# ON THE QUATERNIONIC INVOLUTE-EVOLUTE CURVES

Tülay SOYFİDAN, Mehmet Ali GÜNGÖR

**Abstract:** In this study, after introducing algebraic properties of real quaternions some characterizations of quaternionic involute-evolute curves in $Q$ are obtained. And some results and theorems for quaternionic $w-curves$ are given. Moreover, it is seen that the spatial quaternionic curves in $R^3$ associated with quaternionic involute-evolute curves in $\mathbb{R}^4$ aren't quaternionic involute-evolute curves. Lastly, we illustrate some examples and draw their figures with Mathematica Programme.

## 1. INTRODUCTION

The quaternion was introduced by Irish mathematician Sir William R. Hamilton who discovered that the appropriate generalization is one in which the real axis is left unchanged whereas the vector (imaginary) axis is supplemented by adding two further vector axis in 1843, [5]. Until the middle of the 20th century, the practical use of quaternions was minimal in comparison with other methods. But, currently, this situation has changed. Today, quaternions play a significant role in several areas. In particular for calculations involving three-dimensional rotations such as in three dimensional computer graphics and computer vision. They can be used alongside other methods, such as Euler angles and matrices, or as an alternative to them depending on the application. They can be also used in mechanics, for example, quaternionic formulation of equation of motion in the theory of relativity.

As a set, the quaternions $Q$ are coincide with $R^4$, a four-dimensional vector space over the real numbers. The Serret-Frenet formulae for a quaternionic curves in $\mathbb{R}^3$ and $\mathbb{R}^4$ are introduced by K. Bharathi and M. Nagaraj, [2]. Then, lots of studies have been published by using this studies. One of them is Karadağ and Sivridağ's study which they gave many characterizations for quaternionic inclined curves in $\mathbb{R}^4$, [4]. Another is Gök et al.'s study. They defined a new kind of slant helix in $\mathbb{R}^4$, which they called $\boldsymbol{B}_2-$slant helix and they gave some characterizations of this slant helix in $\mathbb{R}^4$, [10].

The involute-evolute curves in $\mathbb{R}^3$ are well known in elementary differential geometry, see for example, [12]. Moreover, it is well known that if a curve is differentiable at the each point of an open interval then a set of mutually orthogonal unit vectors can be constructed. These vectors are called tangent, normal and binormal unit vectors, often called $\boldsymbol{t},\boldsymbol{n}$ and $\boldsymbol{b}$ or collectively the Serret-Frenet frame. The Serret-Frenet formulas are defined as follows

$$\boldsymbol{t}' = \kappa\boldsymbol{n},\ \boldsymbol{n}' = -\kappa\boldsymbol{t} + \tau\boldsymbol{b},\ \boldsymbol{b}' = -\tau\boldsymbol{n}$$

where " ' " is the derivative with respect to arc-length, $\kappa$ and $\tau$ are the curvatures. The set, whose elements are frame vectors and curvatures of a curve, is called Serret-Frenet apparatus of the curves.

A pair of curves are said to be involute- evolute curves if there exists a one to one correspondence between their points such that one's tangent and the other's principal normal are linear dependent at their corresponding points.

A curve is called $w-$curve if the all curvatures of the curve are constants. In $\mathbb{R}^4$, Özyılmaz and Yılmaz calculated Serret-Frenet apparatus of the involute curve using apparatus of the evolute curve which is $w-$ curve.

The main purpose of this paper is to obtain some characterizations of quaternionic involute-evolute curves in quaternionic space $Q$. To do this, since it is a trivial task to write out the Serret-Frenet formulae of the curve in $\mathbb{R}^4$ using quaternions, firstly it is established some characterizations of quaternionic involute-evolute curves in $\mathbb{R}^4$.

Moreover, it is seen that the spatial quaternionic curves in $\mathbb{R}^3$ associated with quaternionic involute-evolute curves in $\mathbb{R}^4$ aren't quaternionic involute-evolute curves. Lastly, we illustrate some examples and draw their figures with Mathematica Programme.

## 2. PRELIMINARIES

To meet the requirements in this section, the basic elements of the theory of quaternions in the Euclidean space are briefly presented. A more complete elementary treatment can be found in [9].

A real quaternion is defined with $q = d + a\boldsymbol{e}_1 + b\boldsymbol{e}_2 + c\boldsymbol{e}_3$ (or $q = S_q + \boldsymbol{V}_q$ where the symbols $S_q = d$ and $\boldsymbol{V}_q = a\boldsymbol{e}_1 + b\boldsymbol{e}_2 + c\boldsymbol{e}_3$ denote scalar and vector part of $q$, respectively) such that

$$i)\quad \boldsymbol{e}_i \times \boldsymbol{e}_i = -\boldsymbol{e}_4 \qquad (\boldsymbol{e}_4 = +1, \quad 1 \le i \le 3)$$

$$ii)\quad \boldsymbol{e}_i \times \boldsymbol{e}_j = \boldsymbol{e}_k = -\boldsymbol{e}_j \times \boldsymbol{e}_i \qquad (1 \le i, j \le 3)$$

where $(ijk)$ is an even permutation of (123) in the Euclidean space $\mathbb{R}^4$. For every $p, q \in Q$, using these basic products we can now expand the product of two quaternion as $p \times q = S_p S_q - \langle \boldsymbol{V}_p, \boldsymbol{V}_q \rangle + S_p \boldsymbol{V}_q + S_q \boldsymbol{V}_p + \boldsymbol{V}_p \wedge \boldsymbol{V}_q$, where we have used the usual inner and cross products in Euclidean space $\mathbb{R}^3$, [1]. The conjugate of the quaternion $q$ is denoted by $\hat{q}$ and defined $q = S_q - \boldsymbol{V}_q$. Thus, we define symmetric, real valued, non-degenerate, bilinear form $h$ as follows:

$$h : Q \times Q \to R,$$

$$(p, q) \to h(p, q) = \frac{1}{2}(p \times q + q \times p)$$

and it is called the quaternion inner product. The norm of a real quaternion $q$ is $\|q\|^2 = h(q, q) = q \times q = d^2 + a^2 + b^2 + c^2$. If $\|q\| = 1$, then $q$ is called a unit quaternion. It is known that the groups of unit real quaternions and unitary matrices $SU(2)$ are isomorphic. Thus, spherical concepts in $S^3$ such as meridians of longitude and parallels of latitude are explained with assistance elements of $SU(2)$. Furthermore, the elements of $SU(3)$ can match with each element of $S^3$, [11].

The $3-$sphere $S^3 \subset Q$ in quaternionic calculus is like the unit circle $S^1 \subset C$ in complex calculus. In fact, $S^3 = \{ q \in Q \mid \|q\| = 1 \}$ constitutes a group under quaternionic multiplication. $q$ is called a spatial quaternion whenever $q + q = 0$, [2]. Moreover, quaternion product of two spatial quaternions is $p \times q = -\langle p, q \rangle + p \wedge q.$ $q$ is a temporal quaternion whenever $q - q = 0$. Any general $q$ can be written as $q = \frac{1}{2}(q + q) + \frac{1}{2}(q - q),$ [1].

## 3. SOME CHARACTERIZATIONS OF QUATERNIONIC INVOLUTE-EVOLUTE CURVES

The $4-$dimensional Euclidean space $\mathbb{R}^4$ is identified with the space of unit quaternions which is denoted by $\mathbb{Q}$. Let

$$\xi : I \subset R \to Q, \qquad s \to \xi(s) = \sum_{i=1}^{4} \xi_i(s) \boldsymbol{e}_i, \qquad (1 \le i \le 4), \qquad \boldsymbol{e}_4 = 1$$

be a smooth curve defined over the interval $I = [0,1]$. Let the parameter $s$ be chosen such that the tangent $\boldsymbol{T} = \xi'(s)$ has unit magnitude.

**THEOREM 3.1.** *Let $\{\boldsymbol{T}, \boldsymbol{N}, \boldsymbol{B}, \boldsymbol{E}\}$ be the Serret-Frenet frame of the quaternionic curve $\xi = \xi(s)$ in $\mathbb{Q}$. Then, the Serret-Frenet equations are*

$$\begin{aligned} &\boldsymbol{T}'(s)=\kappa(s)\mathbf{N}(s), && \boldsymbol{N}(s)=\boldsymbol{t}(s)\times\boldsymbol{T}(s)\\ &\boldsymbol{N}'(s)=-\kappa(s)\boldsymbol{T}(s)+k(s)\boldsymbol{B}(s), && \boldsymbol{B}(s)=\boldsymbol{n}(s)\times\boldsymbol{T}(s)\\ &\boldsymbol{B}'(s)=-k(s)\boldsymbol{N}(s)+\left(r(s)-\kappa(s)\right)\boldsymbol{E}(s), &&\\ &\boldsymbol{E}'(s)=-\left(r(s)-\kappa(s)\right)\boldsymbol{B}(s), && \boldsymbol{E}(s)=\boldsymbol{b}(s)\times\boldsymbol{T}(s) \end{aligned} \tag{3.1}$$

*where* $\kappa(s)=\left\|\boldsymbol{T}'(s)\right\|$, [2].

It is obtained the Serret-Frenet apparatus for the curve $\xi=\xi(s)$ by using of the Serret-Frenet apparatus for a curve $\alpha=\alpha(s)$ in $\mathbb{R}^3$ where $\alpha$ is a spatial quaternionic curve associated with $\xi$ quaternionic curve and $\{\boldsymbol{t},\boldsymbol{n},\boldsymbol{b}\}$ is the Frenet frame of the spatial quaternionic curve $\alpha$ in [2]. Moreover, there are relationships between curvatures of the curves $\xi$ and $\alpha$ in [2]. These relations can be explained that the torsion of $\xi$, $k$, is the principal curvature of the curve $\alpha$. Also, the bitorsion of $\xi$ is $(r-\kappa)$, where $r$ is the third curvature of $\alpha$ and $\kappa$ is the principal curvature of $\xi$.

Moreover, the Serret-Frenet apparatus of the quaternionic curve $\xi$, are given by

$$\begin{aligned} &\boldsymbol{T}(s)=\frac{1}{\|\xi'(s)\|}\xi'(s),\\ &\boldsymbol{N}(s)=\frac{\|\xi'(s)\|^2\xi''(s)-h\left(\xi'(s),\xi''(s)\right)\xi'(s)}{\left\|\|\xi'(s)\|^2\xi''(s)-h\left(\xi'(s),\xi''(s)\right)\xi'(s)\right\|},\\ &\boldsymbol{B}(s)=\eta\boldsymbol{E}(s)\wedge\boldsymbol{T}(s)\wedge\boldsymbol{N}(s),\\ &\boldsymbol{E}(s)=\eta\frac{\boldsymbol{T}(s)\wedge\boldsymbol{N}(s)\wedge\xi'''(s)}{\|\boldsymbol{T}(s)\wedge\boldsymbol{N}(s)\wedge\xi'''(s)\|},\qquad(\eta=\pm1),\\ &\kappa(s)=\frac{\left\|\|\xi'(s)\|^2\xi''(s)-h\left(\xi'(s),\xi''(s)\right)\xi'(s)\right\|}{\|\xi'(s)\|^4},\\ &k(s)=\frac{\|\boldsymbol{T}(s)\wedge\boldsymbol{N}(s)\wedge\xi'''(s)\|\|\xi'(s)\|}{\left\|\|\xi'(s)\|^2\xi''(s)-h\left(\xi'(s),\xi''(s)\right)\xi'(s)\right\|},\\ &\left(r-\kappa\right)(s)=\frac{h\left(\xi^{(iv)}(s),\boldsymbol{E}(s)\right)}{\|\boldsymbol{T}(s)\wedge\boldsymbol{N}(s)\wedge\xi'''(s)\|\|\xi'(s)\|}. \end{aligned} \tag{3.2}$$

**Definition 3.2.** Let $\phi,\xi:I\subset R\to Q$ be any regular real quaternionic curves with parameter $s^*$ and $s$, respectively. Moreover, $\{\boldsymbol{T}_\phi,\boldsymbol{N}_\phi,\boldsymbol{B}_\phi,\boldsymbol{E}_\phi\}$ and $\{\boldsymbol{T}_\xi,\boldsymbol{N}_\xi,\boldsymbol{B}_\xi,\boldsymbol{E}_\xi\}$ denote the Serret-Frenet frame of the curves $\phi$ and $\xi$, respectively. If $h(\boldsymbol{T}_\phi(s^*),\boldsymbol{T}_\xi(s))=0,$ then, we call curves $\{\phi,\xi\}$ as real quaternionic involute-evolute curves in $\mathbb{Q}$.

**THEOREM 3.3.** *Let* $\xi:I\to Q$ *be a regular unit speed real quaternionic curve. If any real quaternionic curve* $\phi:I\to Q$ *is a quaternionic involute of the curve* $\xi$*, then we have* $d(\xi(s),\phi(s^*))=|c-s|$ *where* $c$ *is real number.*

*Proof.* Let $\phi,\xi:I\to Q$ be a unit speed real quaternionic involute-evolute curves. From definition of quaternionic involute-evolute curves, we write

$$\phi(s^*)=\xi(s)+\lambda(s)\boldsymbol{T}_\xi(s). \tag{3.3}$$

Then differentiating the equation (3.3) with respect to $s$, we obtain that

$$\frac{d\phi}{ds^*}\frac{ds^*}{ds}=\left(1+\lambda'(s)\right)\boldsymbol{T}_\xi(s)+\lambda(s)\kappa_\xi(s)\boldsymbol{N}_\xi(s).$$

By considering the last equation and using quaternionic inner product with tangent vector $\boldsymbol{T}_\xi(s)$, we get

$$0 = h\left(\frac{d\phi}{ds^*}\frac{ds^*}{ds}, \boldsymbol{T}_\xi(s)\right) = \left(1+\lambda'(s)\right) + \lambda(s)\kappa_\xi(s)\left[\frac{1}{2}\left[\boldsymbol{t}(s)\times h\left(\boldsymbol{T}_\xi(s), \boldsymbol{T}_\xi(s)\right) + \left(\boldsymbol{T}_\xi(s)\times \boldsymbol{T}_\xi(s)\right)\times \hat{\boldsymbol{t}}(s)\right]\right]$$

$$= \left(1+\lambda'(s)\right) + \lambda(s)\kappa_\xi(s)\left[\frac{1}{2}\left[\boldsymbol{t}(s) + h\left(\boldsymbol{T}_\xi(s), \boldsymbol{T}_\xi(s)\right)\times \hat{\boldsymbol{t}}(s)\right]\right] = \left(1+\lambda'(s)\right).$$

So, it can be written

$$\lambda(s) = c - s. \tag{3.4}$$

From the equations (3.3) and (3.4), we have

$$d\left(\xi(s), \phi(s^*)\right) = |\lambda| = |c-s|. \tag{3.5}$$

Therefore, from the equations (3.3) and (3.4) we can reach

$$\phi(s^*) = \xi(s) + (c-s)\boldsymbol{T}_\xi(s). \tag{3.6}$$

This completed the proof.

The following theorem provide some characterizations of quaternionic involute-evolute curves.

**THEOREM 3.4.** *Let $\xi : I \to Q$ be a regular unit speed real quaternionic curve and any real quaternionic curve $\phi : I \to Q$ be the involute of $\xi$. The Serret-Frenet apparatus of quaternionic curve $\phi$ can be formed by apparatus of $\xi$.*

*Proof.* Let $\xi = \xi(s)$ be a unit speed real quaternionic curve. Without loss of generality, suppose that $\phi = \phi(s)$ is the involute of $\xi$. By using the equation (3.1), if we calculate the derivatives of the equation (3.6) four times with respect to $s$, we get

$$\frac{d\phi}{ds^*}\frac{ds^*}{ds} = (c-s)\kappa_\xi(s)\boldsymbol{N}_\xi(s), \tag{3.7}$$

$$\begin{aligned}\frac{d^2\phi}{ds^{*2}}\frac{ds^{*2}}{ds^2} &= -(c-s)\kappa_\xi^{\,2}(s)\boldsymbol{T}_\xi(s)\\ &\quad + \left[(c-s)\kappa_\xi'(s) - \kappa_\xi(s)\right]\boldsymbol{N}_\xi(s) + (c-s)\kappa_\xi(s)k(s)\boldsymbol{B}_\xi(s),\end{aligned} \tag{3.8}$$

$$\begin{aligned}\frac{d^3\phi}{ds^{*3}}\frac{ds^{*3}}{ds^3} &= \left[2\kappa_\xi^{\,2}(s) - 3(c-s)\kappa_\xi(s)\kappa_\xi'(s)\right]\boldsymbol{T}_\xi(s)\\ &\quad + \left[-(c-s)\kappa_\xi(s)\left[\kappa_\xi^{\,2}(s) + k^2(s)\right] - 2\kappa_\xi'(s) + (c-s)\kappa_\xi''(s)\right]\boldsymbol{N}_\xi(s)\\ &\quad + \left[2(c-s)\kappa_\xi'(s)k(s) - 2\kappa_\xi(s)k(s) + (c-s)\kappa_\xi(s)k'(s)\right]\boldsymbol{B}_\xi(s)\\ &\quad + (c-s)\kappa_\xi(s)k(s)\left(r(s) - \kappa_\xi(s)\right)\boldsymbol{E}_\xi(s),\end{aligned} \tag{3.9}$$

$$\begin{aligned}\frac{d^4\phi}{ds^{*4}}\frac{ds^{*4}}{ds^4} &= \begin{bmatrix} 9\kappa_\xi(s)\kappa_\xi'(s) - 3(c-s)\kappa_\xi'^{2}(s) - 4(c-s)\kappa_\xi(s)\kappa_\xi''(s)\\ +(c-s)\kappa_\xi^{\,2}(s)\left[\kappa_\xi^{\,2}(s) + k^2(s)\right]\end{bmatrix}\boldsymbol{T}_\xi(s)\\ &\quad + \begin{bmatrix} 3\kappa_\xi^{\,3}(s) - 6(c-s)\kappa_\xi^{\,2}(s)\kappa_\xi'(s) + 3\kappa_\xi(s)k^2(s) - 3(c-s)\kappa_\xi'(s)k^2(s)\\ -3(c-s)\kappa_\xi(s)k(s)k'(s) - 3\kappa_\xi''(s) + (c-s)\kappa_\xi'''(s)\end{bmatrix}\boldsymbol{N}_\xi(s)\\ &\quad + \begin{bmatrix} -6\kappa_\xi'(s)k(s) + 3(c-s)\kappa_\xi''(s)k(s) + 3(c-s)\kappa_\xi'(s)k'(s)\\ -3\kappa_\xi(s)k'(s) + (c-s)\kappa_\xi(s)k''(s) - (c-s)\kappa_\xi(s)k(s)\left(r(s) - \kappa_\xi(s)\right)^2\\ -(c-s)\kappa_\xi(s)k(s)\left[\kappa_\xi^{\,2}(s) + k^2(s)\right]\end{bmatrix}\boldsymbol{B}_\xi(s)\\ &\quad + \begin{bmatrix} 3(c-s)\kappa_\xi'(s)k(s)\left(r(s) - \kappa_\xi(s)\right) - 3\kappa_\xi(s)k(s)\left(r(s) - \kappa_\xi(s)\right)\\ +2(c-s)\kappa_\xi(s)k'(s)\left(r(s) - \kappa_\xi(s)\right) + (c-s)\kappa_\xi(s)k(s)\left(r(s) - \kappa_\xi(s)\right)'\end{bmatrix}\boldsymbol{E}_\xi(s).\end{aligned} \tag{3.10}$$

In view of the equation (3.1), if we take quaternionic norm of the equation (3.7), we obtain that

$$\left\|\frac{d\phi}{ds^*}\frac{ds^*}{ds}\right\|^2 = (c-s)\kappa_\xi(s)\boldsymbol{N}_\xi(s)\times\left((c-s)\kappa_\xi(s)\boldsymbol{N}_\xi(s)\right) = (c-s)^2\kappa_\xi^{\,2}(s)\left(\left(\boldsymbol{t}(s)\times\boldsymbol{T}_\xi(s)\right)\times\left(\boldsymbol{t}(s)\times\boldsymbol{T}_\xi(s)\right)\right)$$

$$= (c-s)^2\kappa_\xi^{\,2}(s)\left(\boldsymbol{t}(s)\times\left(\boldsymbol{T}_\xi(s)\times\boldsymbol{T}_\xi(s)\right)\times\hat{\boldsymbol{t}}(s)\right) = (c-s)^2\kappa_\xi^{\,2}(s)\left(\boldsymbol{t}(s)\times\hat{\boldsymbol{t}}(s)\right) = (c-s)^2\kappa_\xi^{\,2}(s)\left[-\left(-\langle\boldsymbol{t}(s), \boldsymbol{t}(s)\rangle + \boldsymbol{t}(s)\wedge\boldsymbol{t}(s)\ \right)\right].$$

Since $\langle \boldsymbol{t}, \boldsymbol{t}\rangle = 1$, $\boldsymbol{t}\wedge\boldsymbol{t} = 0$ and $\hat{\boldsymbol{t}} = -\boldsymbol{t}$ , we have

$$\left\|\frac{d\phi}{ds^*}\frac{ds^*}{ds}\right\| = \kappa_\xi(s)\left|(c-s)\right| \tag{3.11}$$

and by using the equations (3.2), (3.7) and (3.11) we obtain that

$$\boldsymbol{T}_{\phi}(s^{*}) = \frac{\dfrac{ds^{*}}{ds}(c-s)\kappa_{\xi}(s)}{\left|\dfrac{ds^{*}}{ds}\right|\kappa_{\xi}(s)\left|(c-s)\right|}\boldsymbol{N}_{\xi}(s).$$

So, it is seen that $\boldsymbol{T}_{\phi}(s^{*}) = \pm\boldsymbol{N}_{\xi}(s).$ Without loss of generality let us choose

$$\boldsymbol{T}_{\phi}(s^{*}) = \boldsymbol{N}_{\xi}(s) \tag{3.12}$$

and let us calculate $\boldsymbol{N}_{\phi}(s^{*})$. From the equations (3.7) and (3.8), we obtain that

$$h\left(\frac{d\phi}{ds^{*}}, \frac{d^{2}\phi}{ds^{*2}}\right) = \frac{ds^{3}}{ds^{*3}}\begin{bmatrix} -(c-s)^{2}\kappa_{\xi}^{3}(s)\dfrac{1}{2}\left(\boldsymbol{N}_{\xi}(s)\times\hat{\boldsymbol{T}}_{\xi}(s)+\boldsymbol{T}_{\xi}(s)\times\boldsymbol{N}_{\xi}(s)\right) \\ +(c-s)^{2}\kappa_{\xi}^{2}(s)k(s)\dfrac{1}{2}\left(\boldsymbol{N}_{\xi}(s)\times\boldsymbol{B}_{\xi}(s)+\boldsymbol{B}_{\xi}(s)\times\boldsymbol{N}_{\xi}(s)\right) \\ +(c-s)\kappa_{\xi}(s)\left[(c-s)\kappa_{\xi}'(s)-\kappa_{\xi}(s)\right]\dfrac{1}{2}\left(\boldsymbol{N}_{\xi}(s)\times\boldsymbol{N}_{\xi}(s)+\boldsymbol{N}_{\xi}(s)\times\boldsymbol{N}_{\xi}(s)\right) \end{bmatrix}$$

$$= \frac{ds^{3}}{ds^{*3}}\left[(c-s)^{2}\kappa_{\xi}(s)\kappa_{\xi}'(s)-(c-s)\kappa_{\xi}^{2}(s)\right].$$

Thus by using the equations (3.2), (3.7), (3.8), (3.11) and the last equation we obtain that

$$\boldsymbol{N}_{\phi}(s^{*}) = \frac{-\kappa_{\xi}(s)}{\sqrt{\kappa_{\xi}^{2}(s)+k^{2}(s)}}\boldsymbol{T}_{\xi}(s)+\frac{k(s)}{\sqrt{\kappa_{\xi}^{2}(s)+k^{2}(s)}}\boldsymbol{B}_{\xi}(s). \tag{3.13}$$

By considering the equations (3.9), (3.12) and (3.13) we get

$$\boldsymbol{T}_{\phi}(s^{*})\wedge\boldsymbol{N}_{\phi}(s^{*})\wedge\frac{d^{3}\phi}{ds^{*3}}\frac{ds^{*3}}{ds^{3}} = \frac{\kappa_{\xi}(s)(c-s)}{\sqrt{\kappa_{\xi}^{2}(s)+k^{2}(s)}}\begin{bmatrix} k^{2}(s)\left(r(s)-\kappa_{\xi}(s)\right)\boldsymbol{T}_{\xi}(s) \\ +\kappa_{\xi}(s)k(s)\left(r(s)-\kappa_{\xi}(s)\right)\boldsymbol{B}_{\xi}(s) \\ +\left(\kappa_{\xi}'(s)k(s)-\kappa_{\xi}(s)k'(s)\right)\boldsymbol{E}_{\xi}(s) \end{bmatrix}. \tag{3.14}$$

Therefore, by taking quaternionic norm of the equation (3.14), we obtain that

$$\left\|\boldsymbol{T}_{\phi}(s^{*})\wedge\boldsymbol{N}_{\phi}(s^{*})\wedge\frac{d^{3}\phi}{ds^{*3}}\right\| = \frac{ds^{3}}{ds^{*3}}\frac{\kappa_{\xi}(s)\left|(c-s)\right|}{\sqrt{\kappa_{\xi}^{2}(s)+k^{2}(s)}}\sqrt{\begin{matrix} k^{4}(s)\left(r(s)-\kappa_{\xi}(s)\right)^{2}+\kappa_{\xi}^{2}(s)k^{2}(s)\left(r(s)-\kappa_{\xi}(s)\right)^{2} \\ +\left(\kappa_{\xi}'(s)k(s)-\kappa_{\xi}(s)k'(s)\right)^{2} \end{matrix}}. \tag{3.15}$$

Moreover, using the equations (3.2), (3.14) and (3.15), we find that

$$\boldsymbol{E}_{\phi}(s^{*}) = \eta\frac{k^{2}(s)\left(r(s)-\kappa_{\xi}(s)\right)\boldsymbol{T}_{\xi}(s)+\kappa_{\xi}(s)k(s)\left(r(s)-\kappa_{\xi}(s)\right)\boldsymbol{B}_{\xi}(s)+\left(\kappa_{\xi}'(s)k(s)-\kappa_{\xi}(s)k'(s)\right)\boldsymbol{E}_{\xi}(s)}{\sqrt{k^{4}(s)\left(r(s)-\kappa_{\xi}(s)\right)^{2}+\kappa_{\xi}^{2}(s)k^{2}(s)\left(r(s)-\kappa_{\xi}(s)\right)^{2}+\left(\kappa_{\xi}'(s)k(s)-\kappa_{\xi}(s)k'(s)\right)^{2}}} \tag{3.16}$$

where $\eta = \pm 1$ providing that $\det\left(\boldsymbol{T}_{\xi}(s), \boldsymbol{N}_{\xi}(s), \boldsymbol{B}_{\xi}(s), \boldsymbol{E}_{\xi}(s)\right) = +1$. Similarly, considering the equations (3.2), (3.12), (3.13), (3.16) and necessary arrangements, the binormal vector $\boldsymbol{B}_{\phi}(s^{*})$ is obtained as follows that

$$\boldsymbol{B}_{\phi}(s^{*}) = \eta\frac{\begin{matrix}\left(\kappa_{\xi}'(s)k(s)-\kappa_{\xi}(s)k'(s)\right)\left[k(s)\boldsymbol{T}_{\xi}(s)+\kappa_{\xi}(s)\boldsymbol{B}_{\xi}(s)\right] \\ -k(s)\left(r(s)-\kappa_{\xi}(s)\right)\left(\kappa_{\xi}^{2}(s)+k^{2}(s)\right)\boldsymbol{E}_{\xi}(s)\end{matrix}}{\sqrt{\kappa_{\xi}^{2}(s)+k^{2}(s)}\sqrt{\begin{matrix}k^{4}(s)\left(r(s)-\kappa_{\xi}(s)\right)^{2}+\kappa_{\xi}^{2}(s)k^{2}(s)\left(r(s)-\kappa_{\xi}(s)\right)^{2} \\ +\left(\kappa_{\xi}'(s)k(s)-\kappa_{\xi}(s)k'(s)\right)^{2}\end{matrix}}}. \tag{3.17}$$

Lastly, by using the equation (3.2) if we make necessary arrangement we have the first, second and third curvatures of the quaternionic curve $\phi$ as follows

$$\kappa_{\phi}(s^{*}) = \frac{\sqrt{\kappa_{\xi}^{2}(s)+k^{2}(s)}}{\kappa_{\xi}(s)\left|c-s\right|},$$

$$k^{*}(s^{*}) = \frac{\sqrt{k^{4}(s)\left(r(s)-\kappa_{\xi}(s)\right)^{2}+\kappa_{\xi}^{2}(s)k^{2}(s)\left(r(s)-\kappa_{\xi}(s)\right)^{2}+\left(\kappa_{\xi}'(s)k(s)-\kappa_{\xi}(s)k'(s)\right)^{2}}}{\kappa_{\xi}(s)\left|c-s\right|\left(\kappa_{\xi}^{2}(s)+k^{2}(s)\right)}, \tag{3.18}$$

$$r^*(s^*)-\kappa_\phi(s^*)=\frac{\begin{bmatrix}-(c-s)\kappa_\xi(s)\kappa_\xi''(s)k^2(s)\big(r(s)-\kappa_\xi(s)\big)+2(c-s)\kappa_\xi(s)\kappa_\xi'(s)k(s)k_\xi'(s)\big(r(s)-\kappa_\xi(s)\big)\\ +(c-s)\kappa_\xi^2(s)k(s)k''(s)\big(r(s)-\kappa_\xi(s)\big)-2(c-s)\kappa_\xi^2(s)k'^2(s)\big(r(s)-\kappa_\xi(s)\big)\\ -(c-s)\kappa_\xi^2(s)k^2(s)\big(r(s)-\kappa_\xi(s)\big)^3-(c-s)\kappa_\xi^2(s)k(s)k'(s)\big(r(s)-\kappa_\xi(s)\big)'\end{bmatrix}}{\dfrac{(c-s)^2\kappa_\xi^2(s)}{\sqrt{\kappa_\xi^2(s)+k^2(s)}}\begin{bmatrix}k^4(s)\big(r(s)-\kappa_\xi(s)\big)^2+\kappa_\xi^2(s)k^2(s)\big(r(s)-\kappa_\xi(s)\big)^2\\ +\big(\kappa_\xi'(s)k(s)-\kappa_\xi(s)k'(s)\big)^2\end{bmatrix}}. \tag{3.19}$$

This completed the proof.
From the last theorem we give following corollary.

**COROLLARY 3.5.** Let $\xi : I \to Q$ be a regular unit speed real quaternionic curve. If quaternionic curve $\xi$ is a quaternionic $w-curve$ , then the Serret-Frenet frame of $\xi$ are

$$\begin{aligned}&\boldsymbol{T}_\phi(s^*)=\boldsymbol{N}_\xi(s), && \boldsymbol{N}_\phi(s^*)=\frac{-\kappa_\xi(s)}{\sqrt{\kappa_\xi^2(s)+k^2(s)}}\boldsymbol{T}_\xi(s)+\frac{k(s)}{\sqrt{\kappa_\xi^2(s)+k^2(s)}}\boldsymbol{B}_\xi(s),\\ &\boldsymbol{B}_\phi(s^*)=\eta\boldsymbol{E}_\xi(s), && \boldsymbol{E}_\phi(s^*)=\eta\left[\frac{k(s)}{\sqrt{\kappa_\xi^2(s)+k^2(s)}}\boldsymbol{T}_\xi(s)+\frac{\kappa_\xi(s)}{\sqrt{\kappa_\xi^2(s)+k^2(s)}}\boldsymbol{B}_\xi(s)\right].\end{aligned} \tag{3.20}$$

The following theorem give some characterizations of quaternionic involute-evolute curves.

**THEOREM 3.6.** *Let $\xi,\phi : I \to Q$ be a regular unit speed real quaternionic involute-evolute curves. If the quaternionic curve $\xi$ is a quaternionic $w-curve$, we have*

$$\xi(s)=\phi(s^*)+\rho(s^*)\boldsymbol{N}_\phi(s^*)-\rho(s^*)\frac{k(s)}{\kappa_\xi(s)}\boldsymbol{E}_\phi(s^*),$$

*where $\rho(s^*)$ is the radius of curvature of quaternionic curve $\phi$.*

*Proof.* Suppose that the quaternionic curve $\xi$ is a quaternionic $w-curve$. Since the vector $\xi(s)-\phi(s^*)$ is perpendicular to $\boldsymbol{T}_\phi(s^*)$ as shown in the equation (3.3), we can write

$$\xi(s)=\phi(s^*)+\lambda(s^*)\boldsymbol{N}_\phi(s^*)+\mu(s^*)\boldsymbol{B}_\phi(s^*)+\gamma(s^*)\boldsymbol{E}_\phi(s^*). \tag{3.21}$$

where $\lambda(s^*)$, $\mu(s^*)$ and $\gamma(s^*)$ are arbitrary differentiable functions. If we take the derivative of the equation (3.21) with respect to $s^*$ , then we obtain that

$$\begin{aligned}\boldsymbol{T}_\xi(s)\frac{ds}{ds^*}&=\boldsymbol{T}_\phi(s^*)+\lambda'(s^*)\boldsymbol{N}_\phi(s^*)+\lambda(s^*)\big(\boldsymbol{t}^{*\prime}(s^*)\times\boldsymbol{T}_\phi(s^*)+\boldsymbol{t}^*(s^*)\times\boldsymbol{T}_\phi'(s^*)\big)\\ &\quad+\mu'(s^*)\boldsymbol{B}_\phi(s^*)+\mu(s^*)\big(\boldsymbol{n}^{*\prime}(s^*)\times\boldsymbol{T}_\phi(s^*)+\boldsymbol{n}^*(s^*)\times\boldsymbol{T}_\phi'(s^*)\big)\\ &\quad+\gamma'(s^*)\boldsymbol{E}_\phi(s^*)+\gamma(s^*)\big(\boldsymbol{b}^{*\prime}(s^*)\times\boldsymbol{T}_\phi(s^*)+\boldsymbol{b}^*(s^*)\times\boldsymbol{T}_\phi'(s^*)\big)\\ &=\big(1-\lambda(s^*)\big)\boldsymbol{T}_\phi(s^*)+\big(\lambda'(s^*)+\mu(s^*)k^*(s^*)\big)\boldsymbol{N}_\phi(s^*)\\ &\quad+\Big(\lambda(s^*)k^*(s^*)+\mu'(s^*)-\gamma(s^*)\big(r^*(s^*)-\kappa_\phi(s^*)\big)\Big)\boldsymbol{B}_\phi(s^*)\ +\mu(s^*)\big(r^*(s^*)-\kappa_\phi(s^*)\big)\boldsymbol{E}_\phi(s^*).\end{aligned} \tag{3.22}$$

By using the equations (3.1) and (3.22), we get

$$\begin{aligned}0&=h\left(\boldsymbol{T}_\xi(s)\frac{ds}{ds^*},\boldsymbol{T}_\phi(s^*)\right)\\ &=\big(1-\lambda(s^*)\kappa_\phi(s^*)\big)h\big(\boldsymbol{T}_\phi(s^*),\boldsymbol{T}_\phi(s^*)\big)+\big(\lambda'(s^*)+\mu(s^*)k^*(s^*)\big)h\big(\boldsymbol{N}_\phi(s^*),\boldsymbol{T}_\phi(s^*)\big)\\ &\quad+\Big(\lambda(s^*)k^*(s^*)+\mu'(s^*)-\gamma(s^*)\big(r^*(s^*)-\kappa_\phi(s^*)\big)\Big)h\big(\boldsymbol{B}_\phi(s^*),\boldsymbol{T}_\phi(s^*)\big)\\ &\quad+\mu(s^*)\big(r^*(s^*)-\kappa_\phi(s^*)\big)h\big(\boldsymbol{E}_\phi(s^*),\boldsymbol{T}_\phi(s^*)\big)\\ &=1-\lambda(s^*)\kappa_\phi(s^*).\end{aligned}$$

So, we can write

$$\lambda(s^*)=\frac{1}{\kappa_\phi(s^*)}=\rho(s^*), \tag{3.23}$$

where $\rho(s^*)$ is the radius of curvature of quaternionic curve $\phi.$ From the equation (3.21), we obtain that

$$h\left(\boldsymbol{N}_{\phi}(s^*),\xi(s)-\phi(s^*)\right)=\rho(s^*),$$
$$h\left(\boldsymbol{B}_{\phi}(s^*),\xi(s)-\phi(s^*)\right)=\mu(s^*),$$
$$h\left(\boldsymbol{E}_{\phi}(s^*),\xi(s)-\phi(s^*)\right)=\gamma(s^*)$$

and from Corollary 3.5 and the last equations, we have

$$\rho(s^*)=\frac{-\kappa_{\xi}(s)}{\sqrt{\kappa_{\xi}^{2}(s)+k^{2}(s)}}h\left(\boldsymbol{T}_{\xi}(s),\xi(s)-\phi(s^*)\right)+\frac{k(s)}{\sqrt{\kappa_{\xi}^{2}(s)+k^{2}(s)}}h\left(\boldsymbol{B}_{\xi}(s),\xi(s)-\phi(s^*)\right), \tag{3.24}$$

$$\mu(s^*)=h\left(\boldsymbol{E}_{\xi}(s),\xi(s)-\phi(s^*)\right), \tag{3.25}$$

$$\gamma(s^*)=\frac{k(s)}{\sqrt{\kappa_{\xi}^{2}(s)+k^{2}(s)}}h\left(\boldsymbol{T}_{\xi}(s),\xi(s)-\phi(s^*)\right)+\frac{\kappa_{\xi}(s)}{\sqrt{\kappa_{\xi}^{2}(s)+k^{2}(s)}}h\left(\boldsymbol{B}_{\xi}(s),\xi(s)-\phi(s^*)\right). \tag{3.26}$$

Since $\xi(s)-\phi(s^*)$ is linear dependent with $\boldsymbol{T}_{\xi}$, it is perpendicular to $\boldsymbol{N}_{\xi}(s)$, $\boldsymbol{B}_{\xi}(s)$ and $\boldsymbol{E}_{\xi}(s)$. So, by using the equation (3.25), we obtain that

$$\mu(s^*)=0. \tag{3.27}$$

By considering the equations (3.24) and (3.26) and if we make necessary arrangements, we get

$$\gamma(s^*)=-\frac{\rho(s^*)k(s)}{\kappa_{\xi}(s)}. \tag{3.28}$$

Thus, by the aid of the equations (3.23), (3.27) and (3.28), we can easily write

$$\xi(s)=\phi(s^*)-\rho(s^*)\boldsymbol{N}_{\phi}(s^*)-\frac{\rho(s^*)k(s)}{\kappa_{\xi}(s)}\boldsymbol{E}_{\phi}(s^*). \tag{3.29}$$

This completed the proof.

Let $\xi=\xi(s)$ and $\phi=\phi(s)$ be a unit speed real quaternionic $w-curves$. Thus, we can give the following theorem.

**THEOREM 3.7.** *Let $\xi=\xi(s)$ and $\phi=\phi(s)$ be a regular unit speed real quaternionic $w-curves$ in $\mathbb{Q}$ and $\xi$ is the evolute of the quaternionic curve $\phi$. The Serret-Frenet apparatus of the quaternionic curve $\xi$ can be formed by apparatus of the quaternionic curve $\phi$.*

*Proof.* Suppose that $\xi=\xi(s)$ is the evolute of curve $\phi=\phi(s)$. If we take the derivative of (3.29) with respect to $s^*$, we have

$$\boldsymbol{T}_{\xi}(s)\frac{ds}{ds^*}=\boldsymbol{T}_{\phi}(s^*)-\rho(s^*)\kappa_{\phi}(s^*)\boldsymbol{T}_{\phi}(s^*)+\rho(s^*)k^*(s^*)\boldsymbol{B}_{\phi}(s^*)+\frac{\rho(s^*)k(s)}{\kappa_{\xi}(s)}\left(r^*(s^*)-\kappa_{\phi}(s^*)\right)\boldsymbol{B}_{\phi}(s^*).$$

In view of the equation (3.23), we get

$$\boldsymbol{T}_{\xi}(s)\frac{ds}{ds^*}=\frac{\kappa_{\xi}(s)k^*(s^*)+k(s)\left(r^*(s^*)-\kappa_{\phi}(s^*)\right)}{\kappa_{\xi}(s)\kappa_{\phi}(s^*)}\boldsymbol{B}_{\phi}(s^*).$$

Therefore, it is seen that $\boldsymbol{T}_{\xi}$ is linear dependent with $\boldsymbol{B}_{\phi}$. Without loss of generality, let us choose

$$\boldsymbol{T}_{\xi}(s)=\boldsymbol{B}_{\phi}(s^*). \tag{3.30}$$

As the quaternionic curves $\xi$ and $\phi$ are $w-curves$, $\dfrac{ds}{ds^*}=\dfrac{\kappa_{\xi}(s)k^*(s^*)+k(s)\left(r^*(s^*)-\kappa_{\phi}(s^*)\right)}{\kappa_{\xi}(s)\kappa_{\phi}(s^*)}$ is constant. If we take derivative of the equation (3.30) with respect to $s^*$, we have

$$\frac{d\boldsymbol{T}_{\xi}}{ds}\frac{ds}{ds^*}=\frac{d^2\xi}{ds^2}\frac{ds}{ds^*}=-k^*(s^*)\boldsymbol{N}_{\phi}(s^*)+\left(r^*(s^*)-\kappa_{\phi}(s^*)\right)\boldsymbol{E}_{\phi}(s^*), \tag{3.31}$$

by taking the norm of the last equation, we get

$$\left\|\frac{d^2\xi}{ds^2}\frac{ds}{ds^*}\right\|^2=\begin{bmatrix} k^{*2}(s^*)\left(\boldsymbol{N}_\phi(s^*)\times\boldsymbol{N}_\phi(s^*)\right)-k^*(s^*)\left(r^*(s^*)-\kappa_\phi(s^*)\right)\left(\boldsymbol{N}_\phi(s^*)\times\boldsymbol{E}_\phi(s^*)\right)\\ -k^*(s^*)\left(r^*(s^*)-\kappa_\phi(s^*)\right)\left(\boldsymbol{E}_\phi(s^*)\times\boldsymbol{N}_\phi(s^*)\right)\\ +\left(r^*(s^*)-\kappa_\phi(s^*)\right)^2\left(\boldsymbol{E}_\phi(s^*)\times\boldsymbol{E}_\phi(s^*)\right)\end{bmatrix}$$

$$=\begin{bmatrix} k^{*2}(s^*)h\left(\boldsymbol{N}_\phi(s^*),\boldsymbol{N}_\phi(s^*)\right)\\ -2k^*(s^*)\left(r^*(s^*)-\kappa_\phi(s^*)\right)h\left(\boldsymbol{N}_\phi(s^*),\boldsymbol{E}_\phi(s^*)\right)\\ +\left(r^*(s^*)-\kappa_\phi(s^*)\right)^2 h\left(\boldsymbol{E}_\phi(s^*),\boldsymbol{E}_\phi(s^*)\right)\end{bmatrix}$$

$$=k^{*2}(s^*)+\left(r^*(s^*)-\kappa_\phi(s^*)\right)^2.$$

That is,

$$\left\|\frac{d^2\xi}{ds^2}\frac{ds}{ds^*}\right\|=\sqrt{k^{*2}(s^*)+\left(r^*(s^*)-\kappa_\phi(s^*)\right)^2}. \tag{3.32}$$

By applying the equation (3.2), we get

$$\boldsymbol{N}_\xi(s)=\frac{-k^*(s^*)}{\sqrt{k^{*2}(s^*)+\left(r^*(s^*)-\kappa_\phi(s^*)\right)^2}}\boldsymbol{N}_\phi(s^*)+\frac{\left(r^*(s^*)-\kappa_\phi(s^*)\right)}{\sqrt{k^{*2}(s^*)+\left(r^*(s^*)-\kappa_\phi(s^*)\right)^2}}\boldsymbol{E}_\phi(s^*). \tag{3.33}$$

Similarly, if we take derivative of (3.31) with respect to $s^*$, we have

$$\begin{aligned}\frac{d^3\xi}{ds^3}\frac{ds^2}{ds^{*2}}&=-k^*(s^*)\left(k^*(s^*)\boldsymbol{n}^*(s^*)\times\boldsymbol{T}_\phi(s^*)+\boldsymbol{t}^*(s^*)\times\kappa_\phi(s^*)\boldsymbol{N}_\phi(s^*)\right)\\ &\quad+\left(r^*(s^*)-\kappa_\phi(s^*)\right)\left(-r^*(s^*)\boldsymbol{n}^*(s^*)\times\boldsymbol{T}_\phi(s^*)+\boldsymbol{b}^*(s^*)\times\kappa_\phi(s^*)\boldsymbol{N}_\phi(s^*)\right) \quad (3.34)\\ &=\kappa_\phi(s^*)k^*(s^*)\boldsymbol{T}_\phi(s^*)-\left[k^{*2}(s^*)+\left(r^*(s^*)-\kappa_\phi(s^*)\right)^2\right]\boldsymbol{B}_\phi(s^*).\end{aligned}$$

By using the equations (3.2), (3.30), (3.33) and (3.34), we have

$$\boldsymbol{T}_\xi(s)\wedge\boldsymbol{N}_\xi(s)\wedge\frac{d^3\xi}{ds^3}\frac{ds^2}{ds^{*2}}=\left[-\frac{\kappa_\phi(s^*)k^*(s^*)}{\sqrt{k^{*2}(s^*)+\left(r^*(s^*)-\kappa_\phi(s^*)\right)^2}}\begin{bmatrix}\left(r^*(s^*)-\kappa_\phi(s^*)\right)\boldsymbol{N}_\phi(s^*)\\ +k^*(s^*)\boldsymbol{E}_\phi(s^*)\end{bmatrix}\right]. \tag{3.35}$$

Thus, by taking quaternionic norm of the equation (3.35), we obtain that

$$\left\|\boldsymbol{T}_\xi(s)\wedge\boldsymbol{N}_\xi(s)\wedge\frac{d^3\xi}{ds^3}\right\|=\left|\frac{ds^{*2}}{ds^2}\right|\left|\kappa_\phi(s^*)k^*(s^*)\right|. \tag{3.36}$$

So, from the equations (3.2), (3.35) and (3.36) it can written

$$\boldsymbol{E}_\xi(s)=\eta\left[\frac{\left(r^*(s^*)-\kappa_\phi(s^*)\right)}{\sqrt{k^{*2}(s^*)+\left(r^*(s^*)-\kappa_\phi(s^*)\right)^2}}\boldsymbol{N}_\phi(s^*)+\frac{k^*(s^*)}{\sqrt{k^{*2}(s^*)+\left(r^*(s^*)-\kappa_\phi(s^*)\right)^2}}\boldsymbol{E}_\phi(s^*)\right], \tag{3.37}$$

Similarly, now let us calculate the binormal vector $\boldsymbol{B}_\xi(s)$. If we use the equations (3.2), (3.30), (3.33) and (3.37) and make necessary arrangements we obtain that

$$\boldsymbol{B}_\xi(s)=\eta\sqrt{k^{*2}(s^*)+\left(r^*(s^*)-\kappa_\phi(s^*)\right)^2}\,\boldsymbol{T}_\phi(s^*). \tag{3.38}$$

Lastly, by using the equation (3.2) if we make necessary arrangement we have the first, second and third curvatures of the quaternionic curve $\xi$ as follows

$$\kappa_\xi(s)=\frac{\kappa_\phi(s^*)\sqrt{\left[k^{*2}(s^*)+\left(r^*(s^*)-\kappa_\phi(s^*)\right)^2\right]^3}}{k^*(s^*)\left[k^{*2}(s^*)+\left(r^*(s^*)-\kappa_\phi(s^*)\right)^2+\kappa_\phi(s^*)\left(r^*(s^*)-\kappa_\phi(s^*)\right)\right]},\quad k(s)=\frac{\kappa_\phi{}^2(s^*)\sqrt{k^{*2}(s^*)+\left(r^*(s^*)-\kappa_\phi(s^*)\right)^2}}{k^{*2}(s^*)+\left(r^*(s^*)-\kappa_\phi(s^*)\right)^2+\kappa_\phi(s^*)\left(r^*(s^*)-\kappa_\phi(s^*)\right)}.$$

$$\left(r(s)-\kappa_\xi(s)\right)=\frac{\kappa_\phi{}^2(s^*)\left(r^*(s^*)-\kappa_\phi(s^*)\right)\sqrt{k^{*2}(s^*)+\left(r^*(s^*)-\kappa_\phi(s^*)\right)^2}}{k^*(s^*)\left(k^{*2}(s^*)+\left(r^*(s^*)-\kappa_\phi(s^*)\right)^2+\kappa_\phi(s^*)\left(r^*(s^*)-\kappa_\phi(s^*)\right)\right)}.$$

This calculated the proof.

**THEOREM 4.8.** *$\phi,\xi : I \subset R \to Q$ be regular quaternionic curves with parameter $s^*$ and $s$, respectively. If $(\phi,\xi)$ are quaternionic involute-evolute curves, then the spatial quaternionic curves $(\beta,\alpha)$, associated with $\phi$ and $\xi$, respectively, aren't the spatial quaternionic involute-evolute curves.*

*Proof:* Let $(\phi,\xi)$ be quaternionic involute-evolute curves with parameter $s^*$ and $s$, respectively. The Frenet apparatus of $(\beta,\alpha)$, associated with $\phi$ and $\xi$, are $\{\boldsymbol{t}^*,\boldsymbol{n}^*,\boldsymbol{b}^*,k^*,r^*\}$ and $\{\boldsymbol{t},\boldsymbol{n},\boldsymbol{b},k,r\}$, respectively. So, from the equation (3.1) we write $\boldsymbol{t}^* = \boldsymbol{N}_\phi \times \boldsymbol{T}_\phi$ for the quaternionic curve $\phi = \phi(s^*)$. In there, if we use the equations (3,12) and (3.13), we can write $\boldsymbol{t}^* = (x\boldsymbol{T}_\xi + y\boldsymbol{B}_\xi)\times \boldsymbol{N}_\xi$ where $x = \frac{-\kappa_\xi}{\sqrt{\kappa_\xi{}^2 + k^2}}$, $y = \frac{k}{\sqrt{\kappa_\xi{}^2 + k^2}}$. From the last equation, we can write $\boldsymbol{t}^* = x\left(\boldsymbol{T}_\xi \times \boldsymbol{N}_\xi\right) + y\left(\boldsymbol{B}_\xi \times \boldsymbol{N}_\xi\right)$. By using the equations in (3.1), we find that

$$\boldsymbol{t}^* = x\left(\boldsymbol{T}_\xi \times \boldsymbol{T}_\xi \times \hat{\boldsymbol{t}}\right) + y\left(\boldsymbol{B}_\xi \times \boldsymbol{T}_\xi \times \hat{\boldsymbol{t}}\right) = x\hat{\boldsymbol{t}} + y\left(\boldsymbol{n}\times\hat{\boldsymbol{t}}\right) = -x\boldsymbol{t} + y\boldsymbol{b}.$$

From the last equation, we see that $\boldsymbol{t}^*$ is perpendicular with $\boldsymbol{n}$. However, if $\boldsymbol{t}^*$ is linear dependent with $\boldsymbol{n}$, then $(\beta,\alpha)$ are spatial quaternionic involute-evolute curves. So, since $\boldsymbol{t}^*$ is perpendicular with $\boldsymbol{n}$, $(\beta,\alpha)$ aren't spatial quaternionic involute-evolute curves.

Now, we will give an example for the above theorem.

**Example:** We consider a quaternionic curve with the arc-length parameter $s$, $\xi : I \subset R \to R^4$ as noted below

$$\xi(s) = \left(\cos\left(\frac{s}{2}\right) - \sin\left(\frac{s}{2}\right), \cos\left(\frac{s}{2}\right) + \sin\left(\frac{s}{2}\right), \frac{s}{2}, \frac{s}{2}\right)$$

for all $s \in I$. By considering the equation (3.2) we find that the Frenet apparatus of the quaternionic curve $\xi = \xi(s)$ as follows

$$\boldsymbol{T}_\xi(s) = \frac{1}{2}\left(-\sin\left(\frac{s}{2}\right) - \cos\left(\frac{s}{2}\right), -\sin\left(\frac{s}{2}\right) + \cos\left(\frac{s}{2}\right), 1, 1\right),$$

$$\boldsymbol{N}_\xi(s) = \frac{1}{\sqrt{2}}\left(-\cos\left(\frac{s}{2}\right) + \sin\left(\frac{s}{2}\right), -\cos\left(\frac{s}{2}\right) - \sin\left(\frac{s}{2}\right), 0, 0\right)$$

$$\boldsymbol{B}_\xi(s) = \frac{1}{2}\left(-\cos\left(\frac{s}{2}\right) - \sin\left(\frac{s}{2}\right), \cos\left(\frac{s}{2}\right) - \sin\left(\frac{s}{2}\right), -1, 1\right),$$

$$\boldsymbol{E}_\xi(s) = \frac{1}{\sqrt{2}}(0,0,1,-1).$$

The curvatures functions of the quaternionic curve $\xi = \xi(s)$ are as follows

$$\kappa_\xi = \frac{\sqrt{2}}{4}, \quad k = \frac{\sqrt{2}}{4}, \quad r - \kappa_\xi = 0.$$

By using the equation (3.3), if we make necessary arrangement, we can easily find the quaternionic involute curve of the quaternionic curve $\xi = \xi(s)$, as follows

$$\phi(s) = \frac{1}{2}\left((2-c+s)\cos\left(\frac{s}{2}\right) + (-2-c+s)\sin\left(\frac{s}{2}\right), (2+c-s)\cos\left(\frac{s}{2}\right) + (2-c+s)\sin\left(\frac{s}{2}\right), \frac{c}{2}, \frac{c}{2}\right)$$

which $c$ is a constant.

By using the equation (3.1), the spatial quaternionic curve $\alpha = \alpha(s)$ in $R^3$ associated with quaternionic curve $\xi = \xi(s)$ in $\mathbb{R}^4$ is given by $\alpha(s) = \frac{1}{\sqrt{2}}\left(s, 2\cos\left(\frac{s}{2}\right), 2\sin\left(\frac{s}{2}\right)\right)$ where $s$ is the arc-length parameter of $\alpha$ and its curvature functions are $k = \frac{\sqrt{2}}{4}$, $r = \frac{\sqrt{2}}{4}$. Similarly, the spatial quaternionic curve $\beta$ in $\mathbb{R}^3$ associated with quaternionic curve $\phi = \phi(s)$ in $\mathbb{R}^4$ is given by $\beta(s) = (s+c, c, c)$ which $c$ is a constant.

Now, if we calculate the quaternionic inner product $h\left(\boldsymbol{t}\times\boldsymbol{t}^*\right)$, then we obtain that $h\left(\boldsymbol{t}\times\boldsymbol{t}^*\right)=\frac{1}{\sqrt{2}}\neq 0.$

So, we can easily see that the spatial quaternionic curves $(\beta,\alpha)$ aren't spatial quaternionic involute-evolute curves.

The pictures of the some projections of the quaternionic curve $\xi=\xi(s)$, the quaternionic involute curves of $\xi$ and their associated spatial quaternionic curves are as follows

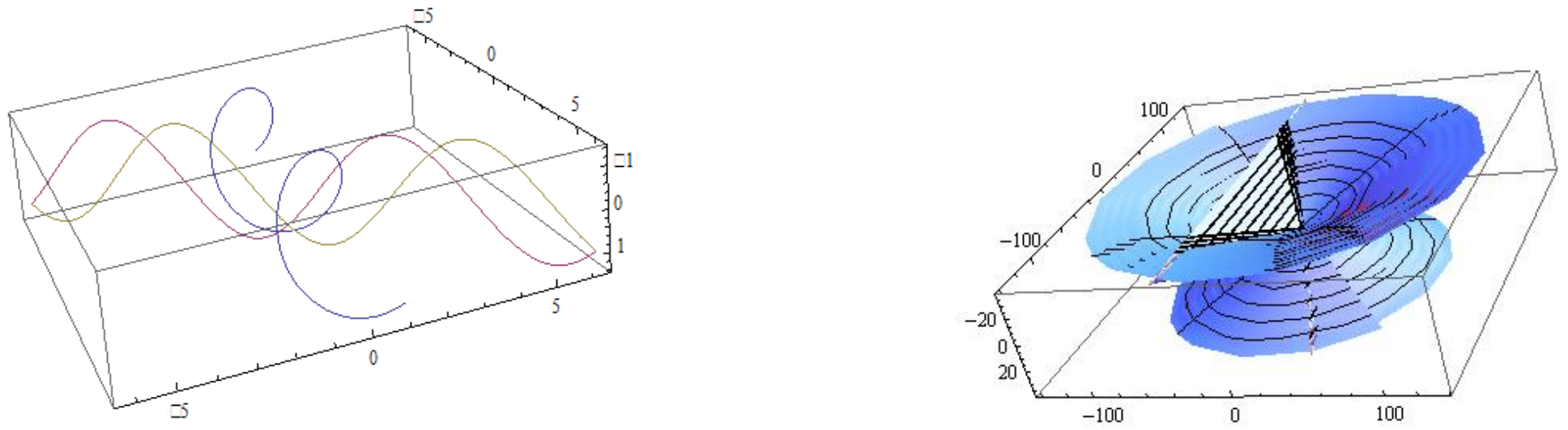

Figure 3.1. Some projections of the quaternionic curve $\xi=\xi(s)$ (on the left) and the quaternionic involutes of $\xi$ (on the right).

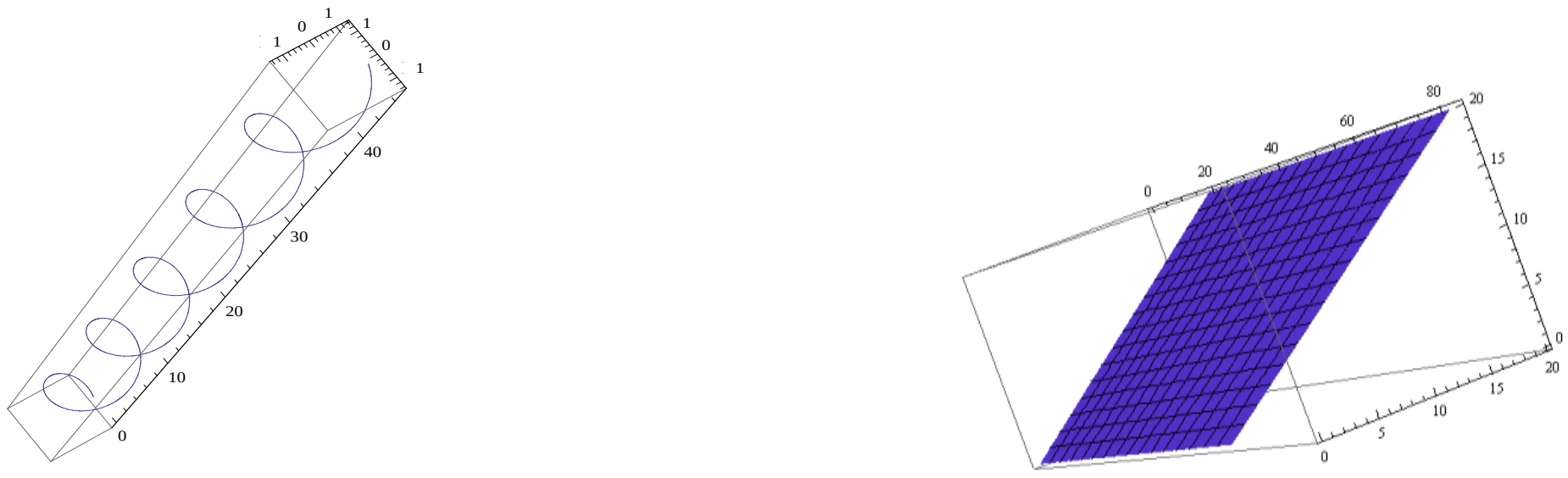

Figure 3.2. The spatial quaternionic curve associated with the quaternionic curve $\xi=\xi(s)$ (on the left) and the spatial quaternionic curves associated with the quaternionic involutes of $\xi$ (on the right).

*Department of Mathematics, Sakarya University, Sakarya, Turkey*
*E-mail: agungor@sakarya.edu.tr*

*Department of Mathematics, Sakarya University, Sakarya, Turkey*
*E-mail: tsoyfidan@sakarya.edu.tr*